# Un énoncé et un texte inaugural

Alain Herreman

*Résumé* : Cet article montre que les principales caractéristiques de la thèse de Turing découlent de celles du système d'expression (machines de Turing) qu'il inaugure, et en particulier de sa *conformité*. Les notions d'*énoncés* et *textes inauguraux* sont à partir de là définies et illustrées par d'autres exemples que la thèse et l'article de Turing. D'autres rapports entre les textes inauguraux sont considérés qui permettent d'inscrire l'article de Turing dans une histoire inaugurale.

*Abstract*: This paper shows that the main features of Turing's thesis derived from those of the expression system (Turing machines) it inaugures, and in particular its conformity. The notions of *inaugural statements* and *texts* are defined and illustrated by other examples than Turing thesis and article. Other relationships between the inaugural texts are considered which allow to include Turing's article in an inaugural history.

## Introduction

Selon la thèse de Turing les nombres tenus pour calculables sont exactement les nombres calculables par les machines de Turing. Le caractère singulier de cet énoncé a d'emblée été reconnu et continu de l'être. Je veux montrer ici qu'il s'agit néanmoins d'un type d'énoncé récurrent en histoire des mathématiques et que je propose d'appeler des *énoncés inauguraux*. Leur caractéristiques découlent de celles du système d'expressions qu'ils inaugurent, la plus remarquable d'entre elles étant sans doute : leur conformité[1]. Ce point de vue consiste à reconnaître la *nécessité* générale d'inaugurer les systèmes d'expressions et à reconnaître ensuite les *problèmes* afférents. Cela amène à distinguer aux côtés des *énoncés* inauguraux des *textes* inauguraux, c'est-à-dire des textes qui présentent des caractéristiques communes en raison de la similitude des problèmes qu'ils doivent résoudre.

---

[1] Kleene qui a introduit le terme spécifique de « thèse » pour désigner cet énoncé donnait aussi l'exemple d'une seconde thèse : « THESIS II. *For any given formal system and given predicate P(a), the predicate that P(a) si provable is expressible in the form (Ex)R(a, x) where R is general recursive.* » Kleene 1943, Davis 1965, 276. La notion générale de thèse a été diversement considérée par Kreisel 1967 ; Barwise 1977 ; Shapiro 1981 ; Epstein & Carnielli 1989, chp 25 ; Mendelson 1990 ; Soare 1996 ; Sieg 1997, 173 ; Black 2000. Mendelson, par exemple, fait entrer dans cette catégorie toute définition d'une notion ayant déjà un sens préalable à la définition mathématique. Il donne l'exemple de la notion de limite définie par Weierstrass ou celle de fonction définie par Peano. Barwise (1977, 41) appelle « thèse de Hilbert », suivant en cela une suggestion de Martin Davis, l'affirmation selon laquelle toutes les propositions et toutes les démonstrations mathématiques peuvent être réécrites au premier ordre. Soares (1996, 296-7) qualifie de « thèse de Cauchy-Weierstrass » l'assertion selon laquelle une fonction est intuitivement continue ssi elle satisfait la définition en ε-δ. Il appelle « thèse des courbes » l'affirmation selon laquelle la notion intuitive de longueur d'une courbe continue dans un espace de dimension 2 est donnée par la définition habituelle de la limite de la somme des longueurs des segments d'une approximation de la courbe. Black (2000, 251) appelle quant à lui « thèse de Hilbert » l'assertion selon laquelle les propositions arithmétiques démontrables coïncident avec les propositions arithmétiques vraies. Il appelle « thèse de Tarski » l'assertion selon laquelle la vérité dans tous les modèles coïncide avec la notion intuitive de validité. La thèse de Cook-Karp désigne couramment l'assertion selon laquelle les problèmes effectivement calculables coïncident avec les problèmes calculables en un temps polynomial. Il n'entre pas dans mon propos de discuter ces différentes acceptions et de les comparer avec celle d'énoncé inaugural.

Dans la première partie de l'article j'introduirai la définition des énoncés et des textes inauguraux. Je montrerai que la thèse et l'article de Turing en sont bien des exemples mais qu'il en existe d'autres. Ces exemples doivent contribuer à convaincre de la pertinence de distinguer, comme tel, ce type d'énoncé et de texte. Mais ils visent surtout à montrer que leurs similitudes doivent être rapportées aux caractéristiques communes des *systèmes d'expressions* que ces différents textes inaugurent.

J'examinerai brièvement dans la seconde partie l'histoire de ces systèmes d'expressions qui, même du point de vue de leur inauguration, ne se réduit pas aux *textes inauguraux*. Il s'agira surtout de mettre en évidence le rapport, historique cette fois, que l'article de Turing entretient à d'autres textes inauguraux.

# I -Textes et énoncés inauguraux

## 1 - Un énoncé inaugural

L'énoncé donné par Turing de la thèse de Turing est le suivant :

> "les nombres calculables [suivant la définition qu'il en donne] incluent tous les nombres que l'on aurait naturellement tendance à considérer comme calculables." (Turing 1936, trad. fr 49)

Cet énoncé présente cinq caractéristiques qui vont servir à définir un *énoncé inaugural* :

1. l'énoncé met en jeu deux totalités (dualisme) ;
2. l'une des totalités est tenue pour pré-établie (réalisme),
3. l'autre totalité n'est pas tenue pour pré-établie sous la forme considérée ou dans le rapport considéré à la totalité pré-établie (inauguration) ;
4. l'énoncé affirme que la deuxième totalité, celle qui n'est pas pré-établie, est une représentation conforme de la première (conformité) ;
5. la démonstration de la conformité des deux totalités est impossible (incommensurabilité[2]) ;

La vérification que l'énoncé de la thèse de Turing est un énoncé inaugural contribuera à préciser le sens et le statut de ces conditions. D'autres exemples seront aussi donnés[3].

---

[2] La notion d'incommensurabilité est ici une notion sémiotique qui ne doit pas être confondue avec celle de Kuhn (1962). Il n'est cependant pas exclu, comme cela pourra apparaître progressivement, de rapporter certains aspects d'une partie des phénomènes considérés par Kuhn à des incommensurabilités sémiotiques. Il importerait ici de distinguer les différentes raisons de l'impossibilité d'établir la conformité.

[3] Une vérification complète pour chacun est ici impossible. La similitude, de ce point de vue, des énoncés de Turing, Church, Post ou Kleene ne posant pas de difficulté je me restreindrai à l'énoncé de Turing, mais tout ce qui est dit ici à propos de l'énoncé de Turing s'applique aux autres. De même, mon propos n'est pas

*a) Le dualisme*

La première condition, le *dualisme*, est la référence à *deux* totalités. C'est bien le cas de la thèse de Turing qui fait référence aux nombres calculables et aux machines logiques. Si ces deux totalités n'étaient pas distinguées, chacune avec une existence propre, la thèse n'aura pas lieu d'être. Et c'est bien en tant que *totalités* que les nombres et les machines logiques sont considérés.
Cette caractéristique se retrouve dans l'énoncé qui ouvre *La Géométrie* de Descartes[4] :

> « Tous les problèmes de Géométrie se peuvent facilement réduire à tels termes, qu'il n'est besoin, par après, que de connaître la longueur de quelques lignes droites, pour les construire. » Descartes 1637, 170

Cet énoncé, comme la suite de *La Géométrie* l'atteste fait référence à deux totalités : les problèmes de Géométrie et les équations polynomiales. C'est aussi un énoncé inaugural.

*b) Le réalisme*

La deuxième condition est le *réalisme*. Ce terme n'implique aucun parti pris philosophique ou ontologique de ma part, il renvoie simplement au *statut objectif dans l'énoncé* de l'une des deux totalités et en l'occurrence à son caractère *pré-établi*.
L'énoncé du Turing suppose qu'*il y a* des nombres calculables. A l'évidence, Turing ne considère pas créer le concept de nombres calculables : il entend seulement donner une définition mathématique d'une notion qu'il considère pré-exister. Il est possible de poser de nombreuses questions sur la nature de cette existence antérieure à la définition mathématique, mais ce critère ne retient que le fait objectif que l'une des deux totalités est tenue pour pré-exister à l'énoncé de la thèse. La totalité qui pré-existe et la nature de celle-ci n'est sans doute pas la même par exemple pour Turing, Church, Kleene et Post. La « philosophie mathématique » qui se manifeste dans leurs commentaires à propos de leurs énoncés est dans certains cas très différente, mais ce caractère pré-établi est bien commun et, de manière variée, souligné. Il ressort par exemple de l'adverbe « naturellement » employé par Turing dans son énoncé (« nombres que l'on aurait naturellement tendance à considérer comme calculables »).
Dans l'énoncé inaugural de Descartes, les problèmes géométriques sont aussi considérés comme une totalité pré-établie. Descartes suppose lui aussi qu'*il y a* des problèmes de géométrie auxquels les mathématiciens font face. Il leur appartient de les résoudre, par de les produire à leur guise. Ils peuvent en découvrir, il ne peuvent pas en inventer. A nouveau, la question n'est pas de savoir si tel est bien en réalité le cas, ni même de savoir si Descartes le croit ou pas sincèrement, mais de reconnaître que son énoncé (et au-delà l'ensemble de *La Géométrie*) le suppose : son énoncé inaugural serait insoutenable, voire vide de sens, si les problèmes de géométrie n'étaient, comme les nombres calculables, considérés comme une totalité constituée et pré-existante.
La représentation des propositions mathématiques par les formules de la logique donne aussi lieu à un énoncé inaugural dans les *Principia mathematica* de Whitehead & Russell :

---

ici de commenter pour elle-même chacune de ces cinq caractéristiques et leurs relations mutuelles mais seulement de vérifier que l'énoncé du Turing les satisfait et d'en donner d'autres exemples.
4 Pour une analyse de ce point de vue de *La Géométrie* de Descartes, voir (Herreman, à paraître).

> "La mathématique pure est la classe de toutes les propositions de la forme "p implique q", où p et q sont des propositions contenant une ou plusieurs variables, les mêmes dans les deux propositions, et où ni p ni q ne contiennent d'autres constantes que des constantes logiques. Et les constantes logiques sont toutes ces notions qui peuvent être définies au moyen de l'implication, de la relation d'un terme à une classe dont il est membre, de la notion de tel que, de la notion de relation, et de toutes les autres notions que peut impliquer celle, générale, de proposition de cette forme." *Principa mathematica*, trad. fr. Russell & Roy 1989, 21

Cette représentation des propositions mathématiques contribue à établir que les mathématiques se réduisent à la logique. Ce logicisme est lui-même fondé, chez Russell, sur un «atomisme logique » selon lequel les propositions logiques s'analysent en atomes logiques comme la matière s'analyse en atomes de matière. Un langage logiquement parfait, comme l'est celui des *Principia Mathematica*, est donc tel qu'aux atomes de ce langage, en l'occurrence ses mots, correspond des choses elles-mêmes simples qui *pré-existent* à la représentation introduite :

> « Dans un langage logiquement parfait, les mots d'une proposition correspondraient un à un aux composants du fait correspondant, à l'exception des mots tels que « ou », « non », « si », « alors », qui ont une fonction différente. Dans un langage logiquement parfait, il y aura un mot et un seul pour chaque objet simple, et tout ce qui n'est pas simple sera exprimé au moyen d'une combinaison de mots, d'une combinaison dérivée bien entendu des mots représentant les choses simples qui entrent dans sa composition, à raison d'un pour chaque composant simple. Un tel langage sera complètement analytique et montrera immédiatement la structure du fait affirmé ou nié. Le langage exposé dans les *Principia Mathematica* prétend être un langage de cette espèce. » Russell 1918, trad. fr. Russell & Roy 1989, 356.

La référence répétée à des « faits » manifeste le caractère pré-établi de ce qui est représenté par le système d'expressions inauguré. Russell peut ainsi assimiler les recherches en logique à celles menées en astronomie[5], en zoologie[6], etc., et dénoncer les systèmes créés de manière *ad hoc* (Russell 1961, 99).

### c) L'inauguration

Un énoncé inaugural fait référence à deux totalités (dualité), dont une a un caractère pré-établi (réalisme). La troisième condition, l'inauguration, consiste au contraire à reconnaître que la seconde totalité n'est pas pré-établie (sous la forme considérée ou dans le rapport considéré à la première).
Il s'agit là simplement de constater que la notion de machines logiques est bien produite par Turing. Si Turing ne prétend pas avoir créé la notion de nombres calculables (ce que l'on pourrait pourtant considérer qu'il a fait...), il peut en revanche revendiquer la définition des machines logiques auxquelles son nom a ensuite été donné.
De même, Descartes ne prétend pas avoir créé la notion de problèmes de géométrie (ce que l'on pourrait aussi considérer qu'il a fait...), il prétend bien en revanche inaugurer leur représentation par des équations polynomiales. Et l'opinion que chacun peut avoir sur l'antériorité de Viète ou de Harriot ne change rien au fait que tel est bien le statut objectif de la représentation par des équations polynomiales *dans La Géométrie*. De même, Whitehead & Russell ne revendiquent pas plus l'invention des propositions logiques que celle des atomes qui, pour Russell, en fondent l'analyse. Ils revendiquent seulement l'inauguration de la représentation logique qu'ils en donnent (ce qui ne les empêchent pas de reconnaître leurs dettes à l'égard de Frege ou de Peano). La revendication de la nouveauté de la représentation introduite (inauguration) suppose la non revendication de la nouveauté de l'autre, c'est-à-dire son caractère pré-établi (réalisme).

---

[5] « La discussion des indéfinissables – qui forme la principale partie de la logique philosophique – consiste en un effort pour voir clairement, et pour faire voir aux autres clairement, les entités dont il est question, de façon que l'esprit puisse avoir d'elles une connaissance du même ordre que celle qu'il a du rouge ou du goût de l'ananas. Là où, comme dans le cas présent, les indéfinissables sont obtenus originairement comme le résidu nécessaire d'une analyse, il est souvent plus facile de savoir qu'il doit y avoir de telles entités que de les percevoir réellement ; le processus est analogue à celui qui aboutit à la découverte de Neptune, avec cette différence que le stade final - la recherche avec un télescope mental de l'entité qui a été déduite – est souvent la partie la plus difficile de l'entreprise. » Russell 1937, xv.
[6] « Logic is concerned with the real world just as truly as zoology » Russell 1919, 169.

*d) Conformité*

La quatrième condition est la conformité, c'est-à-dire l'affirmation selon laquelle la totalité inaugurée donne une représentation conforme de la totalité pré-établie.
Cette condition est à nouveau satisfaite par la thèse de Turing qui affirme que les machines logiques offrent une représentation conforme des nombres calculables[7]. Il faut bien en apprécier la portée : il s'agit non seulement de dire qu'à tout nombre correspond une machine logique, celle qui le calcule, mais inversement que toute machine logique définit elle-même un nombre calculable. Il n'y a ici aucun contenu formel (Granger 1980) : toute machine logique est censée définir un nombre calculable. La représentation représente tout ce qui doit être représenté et elle n'introduit rien qui n'existait pas déjà avant son inauguration : elle n'introduit aucun nouveau nombre calculable. Son introduction est transparente. Mais la conformité suppose encore plus que cela : elle signifie que *toutes les propriétés*, ici des nombres calculables, sont restituées par la représentation inaugurée. Elle ne doit pas seulement représenter les termes, elle doit aussi en reproduire *toutes* les *propriétés*.
Cette conformité se retrouve dans la représentation des fonctions par les séries trigonométriques inaugurée par Fourier. Voici d'abord l'une des formulations de son énoncé inaugural:

> "Il résulte de mes recherches sur cet objet que les fonctions arbitraires même discontinues peuvent toujours être représentées par les développements en sinus ou cosinus d'arcs multiples, et que les intégrales qui contiennent ces développements sont précisément aussi générales que celles où entrent les fonctions arbitraires d'arcs multiples." Fourier 1807, cité *in* Hérivel 1980, 56

L'exigence de conformité est aussi explicite :

> « Les intégrales que nous avons obtenues ne sont point seulement des expressions générales qui satisfont aux équations différentielles : elles représentent de la manière la plus distincte l'effet naturel, qui est l'objet de la question. C'est cette condition principale que nous avons eue toujours en vue, et sans laquelle les résultats du calcul ne nous paraîtraient que des transformations inutiles. Lorsque cette condition est remplie, l'intégrale est, à proprement parler, l'équation du phénomène ; elle en exprime clairement le caractère et le progrès, de même que l'équation finie d'une ligne ou d'une surface courbe fait connaître toutes les propriétés de ces figures. (...) En général, on ne pourrait apporter aucun changement dans la forme de nos solutions sans leur faire perdre leur caractère essentiel, qui est de représenter

---
[7] A la lettre, l'énoncé du Turing n'affirme qu'un sens de l'inclusion, c'est que l'autre (que les machines logiques définissent des nombres calculables) est tenue pour évident. Mais c'est bien la correspondance entre machines logiques et nombres calculables qui est affirmée.

les phénomènes. » Fourier 1822 art. 428

Le « Discours préliminaire » par lequel s'ouvre la *Théorie analytique de la chaleur* de Fourier est en grande partie consacré à défendre plus généralement la *conformité* de l'Analyse mathématique avec les propriétés physiques (Fourier 1822, art. 20).

Russell n'a semble-t-il quant à lui jamais cessé de soutenir la conformité de la logique et des mathématiques. En 1937, il écrit par exemple dans la préface de la réédition de ses *Principles of mathematics* :

> « La thèse fondamentale des pages qui suivent, que les mathématiques et la logique sont identiques, est une thèse que je n'ai depuis jamais vu la moindre raison de modifier. » Russell 1937 cité *in* Grattan-Guinness 1984, 76

Cette conformité est sans doute la condition la plus remarquable et celle qui a le plus de conséquences épistémologiques et historiographiques. En particulier, les définitions mathématique introduites (machines de Turing, équations polynomiales, séries trigonométriques, formules logiques, etc.) ne prétendent pas être « meilleures » que l'acception intuitive (plus précise, plus rigoureuse, sans ambiguïté, etc.), dès lors qu'elles revendiquent avant tout leur *coïncidence*. Leur intérêt est singulièrement de coïncider exactement avec elles et leur valeur tient en partie à cette adéquation. Leur intérêt n'est donc pas, en un sens, de bien définir, ou de permettre de mieux distinguer les fonctions calculables des autres (par exemple de décider de cas douteux), mais de définir la même chose autrement et en l'occurrence d'en donner une représentation qui, bien que nouvelle, pourra être reconnue comme étant de nature mathématique. Elles opèrent le passage d'une définition tenue pour constituée et naturelle, à laquelle aucun reproche ne peut être adressé quant à son extension, aucune rectification n'en sera faite, vers une définition mathématique. La conformité revendiquée contribue à neutraliser tout ce qu'une nouvelle définition, et en l'occurrence une nouvelle représentation, peut avoir de contingent et en particulier de subjectif et d'historique.

### *e) L'incommensurabilité*

La cinquième et dernière condition est l'*incommensurabilité*, c'est-à-dire l'impossibilité de démontrer la conformité des deux totalités. Cette condition est une conséquence des précédentes. Les conditions de dualisme et de réalisme rendent généralement impossible une démonstration ne serait-ce que de la correspondance.
Ainsi, Turing ne peut tout simplement pas démontrer qu'à tout nombre calculable correspond une machine logique parce qu'il ne dispose pas d'une représentation uniforme des nombres calculables le permettant. L'intérêt des machines logiques tient pour une part essentielle dans l'introduction d'une telle représentation qui offre des possibilités de démonstrations impossibles avant cela. Turing signale clairement l'impossibilité de donner une justification mathématique satisfaisante de cette conformité :

> « Les arguments que nous pourrons donner doivent, par principe, faire appel à l'intuition, et seront pour cette raison plutôt insatisfaisants, mathématiquement parlant. » Turing 1936, trad. fr. 76

Tarski (1931) inaugure la représentation des ensembles définissables de nombres par les ensembles de nombres élémentairement définissables. En voici l'énoncé inaugural :

> «tout ensemble particulier des nombres, auquel on a affaire dans les mathématiques, est un ensemble définissable, puisque nous n'avons d'autre moyen d'introduire individuellement un ensemble donné dans le domaine des considérations que celui de construire la fonction propositionnelle qui le détermine, et cette construction constitue par elle-même la preuve de la définissabilité de cet ensemble. » Tarski 1931, 220

L'incommensurabilité le conduit à faire des observations similaires à ceux de Turing :

> "Or, la question s'impose si la définition [d'ensembles (élémentairement) définissables de nombres] qui vient d'être construite et dont la rigueur formelle n'éveille aucune objection, est également juste au point de vue matériel; en d'autres mots, saisit-elle en effet le sens courant et intuitivement connu de la notion? Cette question ne contient, bien entendu, aucun problème de nature purement mathématique, mais elle est néanmoins d'une importance capitale pour nos considérations.» Tarski 1931, 229-230

Peu avant, en 1930, Jacques Herbrand faisait des observations semblables à propos de la conformité de la représentation des propositions et des démonstrations mathématiques proposée par Whitehead & Russell :

> « L'énoncé de toutes ces règles est tel, comme on s'en rendra compte plus loin, que toute combinaison de signes appelée proposition est effectivement la traduction d'une proposition du langage ordinaire, et qu'une qu'une suite de propositions qu'il faut considérer comme vraies d'après ces règles correspond à une démonstration mathématique. Mais on peut encore se demander s'il n'y a pas des propositions que le système des signes ne saura traduire, ou des démonstrations qui échapperont à ses règles.
> Il ne faut pas cacher que cela n'est qu'un résultat expérimental (dont la validité ne peut être que confirmée par une dialectique philosophique) ; sa preuve réside, en somme, dans le fait que Russell et Whitehead ont réussi, dans les trois tomes des *Principia mathematica*, à reproduire tous les raisonnements des débuts des mathématiques et de la théorie des ensembles. On peut considérer comme un des faits les plus parfaitement vérifiés dans notre connaissance logique du monde que tout raisonnement que peut actuellement faire un mathématicien raisonnable trouve immédiatement sa traduction dans le système de signes étudié. Il est possible, quelque difficile que cela nous paraisse, que l'on pourra quelque jour faire éclater ces cadres ; on se convaincra peut-être par la compréhension plus précise de ce qu'est cette logique symbolique que, sans doute, même alors, il lui suffira de modifier quelque peu les règles initiales de l'emploi de ses signes pour absorber des nouvelles théories. » Herbrand 1930 *in* Herbrand 1968, 36

Les cinq conditions retenues, dont on vient de vérifier qu'elles étaient satisfaites par la thèse de Turing, empêchent de confondre les énoncés inauguraux avec l'un quelconque des types d'énoncés habituellement distingués, que ce soit les définitions, les théorèmes ou les axiomes, ou de les réduire simplement à des hypothèses (ce qu'ils sont à l'évidence, mais dont il importe de reconnaître la spécificité).
Le problème posé par l'identification de ces énoncés qui échappent aux catégories reconnues ressort particulièrement bien de l'article dans lequel Cantor inaugure la représentation des points d'une droite par les nombres réels (selon la définition ensembliste qu'il en donne) :

> « Mais pour achever de faire reporter le lien que nous observons entre les systèmes des grandeurs numériques définies dans le §1 et la géométrie de la ligne droite, il faut ajouter encore un axiome dont voici le simple énoncé : A chaque grandeur numérique appartient aussi, réciproquement, un point déterminé de la droite, dont la coordonnée est égale à cette grandeur numérique dans le sens exposé dans ce §. » Cantor 1872, 128, trad. fr. 342

Cantor qualifie ici « d'axiome » son énoncé inaugural et observe à son propos :

> « J'appelle ce théorème un axiome, par ce qu'il est dans sa nature de ne pouvoir être démontré d'une façon générale. Ce théorème sert aussi à donner supplémentairement aux grandeurs numériques une certaine objectivité, dont elles sont, toutefois, complètement indépendants»

L'impossibilité d'établir la conformité des représentations de chacune des deux totalités, due à leur hétérogénéité, à l'absence d'une représentation uniforme aussi bien de la première totalité que des *propriétés* à considérer, oblige celui qui veut soutenir un énoncé inaugural à des développements spécifiques qui permettent de caractériser les *textes inauguraux*.

## 2 - Un texte inaugural

### a) Définition d'un texte inaugural

La conformité de la représentation des nombres calculables par des machines logiques n'est pas évidente ; elle doit être soutenue. Cela requiert des arguments d'un type particulier. Il est aussi objectivement impossible de le faire complètement (incommensurabilité). Les moyens d'expression qui le permettraient font défaut ; il s'agit précisément de les introduire. La conjonction de cette nécessité et de cette impossibilité conduit à un texte d'un genre particulier que j'appelle un *texte inaugural.* En reprenant la caractérisation d'un énoncé inaugural on peut caractériser un tel inaugural de la manière suivante[8] :

1. il met en jeu deux totalités (dualisme) ;
2. l'une des totalités est tenue pour pré-établie (réalisme) ;
3. l'autre totalité n'est pas tenue pour pré-établie sous la forme considérée ou dans le rapport considéré à la totalité pré-établie (inauguration) ;
4. une fonction du texte est de soutenir que la deuxième totalité, celle qui n'est pas pré-établie, est une représentation conforme de la première (conformité) ;
5. la démonstration de la conformité des deux totalités est impossible (incommensurabilité) ;

Turing consacre deux sections entières de son article à soutenir cette conformité. Ce sont ces sections, nécessaires à son propos, qui en font un texte inaugural. Il distingue lui-même clairement trois types d'arguments :

> "Je défendrai mon point de vue au moyen de trois types d'arguments :
> (a) en faisant directement appel à l'intuition ;
> (b) en démontrant l'équivalence de deux définitions (au cas où la nouvelle définition aurait un sens intuitif plus évident) ;
> (c) en exhibant certaines grandes classes de nombres calculables" Turing 1936, trad. 77

Je vais reprendre ses distinctions pour montrer la similitude de sa démarche avec celle d'autres textes inauguraux

---

[8] Il n'y a bien sûr aucune raison de supposer qu'un texte inaugural soit entièrement ou exclusivement inaugural.

### *b) Définitions équivalentes*

Le premier type d'argument qu'il invoque consiste à dégager les opérations élémentaires effectuées par une personne qui calcule et à montrer que ce sont celles effectuées par une machine logique. La description donnée suggère que c'est par une telle analyse que Turing est arrivé à sa définition.

Le deuxième type d'argument consiste à établir l'équivalence de la définition proposée avec une autre. La démonstration de l'équivalence avec la définition de Church en est un exemple. Mais Turing ne disposait pas de cette définition au moment d'écrire son article, et cette démonstration a été ajoutée en appendice. Son argument fait en l'occurrence référence à une autre définition qu'il introduit pour l'occasion et qui utilise le calcul fonctionnel de Hilbert. Je reviendrai plus loin sur cette définition. Il donne ainsi une autre définition dont il peut, cette fois, démontrer l'équivalence avec celle qu'il propose.

*La Géométrie* Descartes comprend un deuxième énoncé inaugural aux courbes géométriques[9] (et non plus relatif aux problèmes) :

> « Et en quelque autre façon qu'on imagine la description d'une ligne courbe, pourvû qu'elle soit du nombre de celles que je nomme Géométriques, on pourra toujours trouver une équation pour déterminer tous ses points en cette sorte. » Descartes 1637, 195

C'est un énoncé sur la représentation des courbes géométriques au moyen d'équations polynomiales. Pour soutenir la conformité de cette représentation Descartes introduit aussi une définition des courbes géométriques au moyen d'instruments dont il s'attache ensuite à établir l'équivalence avec celle par des équations.

La *Begriffsschrift* de Frege est aussi un texte inaugural. Son idéographie s'offre comme une représentation conforme des lois logiques valides :

> « De cette manière on parvient à un petit nombre de lois dans lesquelles, si l'on ajoute celles qui sont contenues dans les règles, le contenu de toutes, bien que latent, est inclus.» Frege 1879, §13, trad. Rivenc & de Rouilhan 1992.

Whitehead et Russell proposent eux-mêmes, expressément à la suite de Frege, une représentation conforme des déductions logiques en même temps que des propositions mathématiques :

> « tout système déductif doit contenir dans ses prémisses autant de propriétés de l'implication nécessaires pour légitimer la procédure habituelle de déduction. Dans cette section, nous poserons comme prémisses certaines propositions et nous montrerons qu'elles sont

---

[9] *La Géométrie* de Descartes inaugure en fait *quatre* représentations. C'est à ma connaissance le seul texte qui en contienne autant.

suffisantes pour toutes les formes courantes d'inférence. » *Principia Mathematica*, Section A, 94 (trad. AH)

D'autres représentations en ont encore ensuite été proposées par Hilbert-Bernays (1934) ou Gentzen (1934). Gentzen (1934) s'est attaché à établir l'équivalence de sa représentation avec celle d'Hilbert & Ackermann (1928), inspirée de celle de Whitehead & Russell[10]. Huntington (1935) s'est attaché à établir l'équivalence de l' « implication matérielle » des *Principia Mathematica* et de l' « implication » de Hilbert-Bernays (1934).

A la différence des énoncés inauguraux, toutes ces équivalences sont bien des théorèmes : aucune incommensurabilité n'empêche plus qu'une démonstration en soit donnée. Mais ce ne sont pas seulement des théorèmes, ce sont aussi des arguments pour soutenir que les définitions proposées sont indépendantes des circonstances de leur énonciation. Elles servent à écarter les objections fondées sur les particularités de chacune. Chaque définition transfère ainsi sur les autres ce que sa formulation reprend de la notion intuitive considérée, et voit comme effacé par les autres ce qu'elle peut y introduire d'étranger. Ces équivalences soutiennent une sorte d'innocuité des représentations introduites. Elles contribuent à rendre transparentes les conditions de leur formulation. Toutes les équivalences qu'il est possible de démontrer entre ces diverses définitions renforcent la conviction que les diverses formulations sont conformes à ce qu'elles définissent en établissant qu'elles sont indépendantes de leur « appareillage ». Elles constituent toujours des arguments effectifs forts. Elles apparaissent d'ailleurs comme des moments déterminants dans l'adhésion à la thèse de Turing (Davis 1982 ; Sieg 1997). Elles ne sauraient néanmoins prétendre établir complètement que ces définitions sont indépendantes des caractéristiques communes des langages qu'elles utilisent.

*c) Exemples cruciaux et classes d'exemples*

Le troisième type d'argument distingué par Turing consiste à vérifier que des nombres que chacun tient pour calculables le sont aussi suivant sa définition. On peut distinguer ici deux procédés.

Le premier consiste à choisir des *exemples cruciaux*. Ce sont des nombres calculables particuliers dont la vérification contribue à convaincre qu'il devrait en être de même pour les autres. Turing démontre ainsi que π et *e* sont calculables suivant sa définition.

Le deuxième procédé consiste à établir la conformité de la définition pour des *classes* entières de nombres. Il s'agit cette fois de vérifier la thèse sinon pour tous les nombres calculables, tout du moins pour certaines grandes familles. Turing le fait pour les "From what has been said above, it should be clear that Gentzen's systems of natural deduction are not arbitrary formalizations of first order logic but constitute a significant analysis of the proofs in this logic.
The situation may be compared to the attempts to characterize the notion of computation where e.g. the formalism of μ-recursive functions or even the general recursive functions may be regarded as an extensional characterization of this notion while Turing's analysis is such that one may reasonably assert the thesis that every computation when sufficiently analysed can be broken down in the operations described by Turing" Prawitz 1993, 209 nombres algébriques réels et les nombres de Bessel.

---

10 Le rapprochement entre l'équivalence, et le rôle propre, des différentes définitions de la calculabilité et celle entre les différentes formulations des règles de déduction a déjà été opéré par Prawitz :

L'un ou l'autre de ces deux procédés se retrouve inévitablement dans tous les textes inauguraux. Dans *La Géométrie*, le problème de Pappus sert ainsi de « problème crucial » : montrer que ce problème qui n'avait jamais été résolu peut être mis en équation et résolu est un argument fort en faveur de la possibilité de représenter et de résoudre algébriquement tous les problèmes de géométrie. Ce qui est aussi remarquable en l'occurrence, c'est que ce problème formulé pour un nombre quelconque de droites, est à la fois *un* problème particulier et lui-même une *classe* de problèmes : le résoudre c'est résoudre une infinité de problèmes couvrant *tous* les genres de problèmes distingués par Descartes.

Ces deux procédés se retrouvent aussi pour soutenir dans *La Géométrie* la conformité des courbes géométriques et des équations polynomiales, par exemple en montrant comment retrouver à partir de leurs équations les paramètres à partir desquels Apollonius construit les coniques.

Fourier ne peut pas non plus échapper aux exemples cruciaux pour soutenir la possibilité de représenter n'importe quelle fonction par une série trigonométrique. Il s'attache ainsi à montrer que le développement en série trigonométrique de la fonction logarithme est donné par l'état permanent de la température dans une barre métallique chauffée à une extrémité :

> « Lorsqu'une barre métallique est exposée par son extrémité à l'action constante d'un foyer et que tous ses points ont acquis leur plus haut degré de chaleur, le système des températures fixes correspond exactement à une Table de logarithmes ; les nombres sont les élévations des thermomètres placés aux différents points, et les logarithmes sont les distances de ces points au foyer.» Fourier 1822, art. 20

Fourier soutient la conformité de sa représentation en montrant qu'il est possible de développer en série trigonométrique des exemples cruciaux comme les fonctions constantes, la fonction identité, cos(x), sin(x) etc.

De même, Frege s'attache à montrer qu'il peut écrire avec la conditionnalité et la négation le « ou », le « ou » exclusif, le « et », l'incompatibilité, le « ni ni » et les quatre oppositions logiques d'Aristote. En établissant qu'il peut formuler et démontrer dans son idéographie le principe d'induction complète, Frege non seulement réduit l'arithmétique à la logique mais il donne aussi un exemple crucial pour son idéographie.

Les trois imposants volumes des *Principia Mathematica* sont aussi, comme on a vu Herbrand le souligner, une preuve par l'exemple que leur représentation des démonstrations et des propositions est conforme aux démonstrations et aux propositions mathématiques.

Tarski est aussi obligé de recourir à de tels exemples pour soutenir la conformité de sa représentation des ensembles de nombres définissables. La similitude de ses propos avec ceux de Turing fait à nouveau ressortir la similitude des problèmes auxquels ils sont confrontés :

> "Si nous désirons acquérir la certitude subjective de la justesse matérielle de la déf. 10 et de sa conformité à l'intuition, sans sortir du domaine des considérations strictement mathématiques, nous sommes contraints de recourir à la voie empirique. On constate notamment, en examinant divers ensembles particuliers qui sont élémentairement définis dans le sens intuitif de cette notion, qu'ils appartiennent tous à la famille $D_1$ ; réciproquement, on parvient à construire une définition élémentaire pour tout ensemble particulier de cette famille: on s'aperçoit aisément que la même méthode de raisonnement, tout à fait automatique, est applicable dans tous les cas." Tarski 1931, 229-230

Une autre représentation conforme est celle des nombres entiers par le système de numération décimale de position. Le texte d'al-Khwarizmi ayant été perdu, citons cet énoncé extrait du *Liber Alchorismi* :

"Puisqu'on a donc montré de quelle manière tout nombre peut être représenté grâce à la disposition des neuf figures citées et du dixième, le zéro, et aussi ce que sont les nombres formés d'unités, formés de dizaines, et composés, il faut passer aux chapitres suivants."
*Liber Alchorismi* in Allard 1992, 74

Les textes qui introduisent en Europe ce système de numération (Allard 1992, Folkerts, 1997) sont essentiellement consacrés à donner la présentation des algorithmes de calcul par lesquels les opérations habituelles sur les nombres peuvent être réalisées (addition, soustraction, division, racines carrées, etc.). Ces textes s'attachent à montrer la conformité de cette représentation en montrant qu'elle permet de noter tous les nombres et d'effectuer les opérations effectuées sur eux. Là aussi, il faut bien apprécier le fondement de notre conviction qu'il est possible d'exprimer ainsi *tous* les nombres entiers et que les algorithmes pour l'addition, la soustraction, la multiplication, etc. sont corrects, c'est-à-dire qu'ils donneront bien toujours le bon résultat.

On peut ainsi constater que Turing développe pour soutenir sa thèse des arguments qui se retrouvent dans les autres textes inauguraux. Tous ces textes ont leur intérêt propre (ils n'introduisent pas les mêmes représentations!). Tous correspondent à des moments importants de l'histoire des mathématiques : disposer d'une représentation qui peut prétendre à la conformité est sans doute trop remarquable pour qu'un tel texte inaugural ne soit pas important. La reconnaissance à partir de ces exemples de ces caractéristiques communes peut contribuer à changer la perspective dans laquelle la thèse de Turing est souvent considérée. Ces similitudes peuvent aussi facilement être relativisées, et l'on peut, sans difficultés, chercher à faire ressortir leurs différences pour récupérer la singularité, réelle, de chacun d'eux. Mais c'est alors en manquer la raison. Si ces textes présentent de telles similitudes c'est parce que leurs auteurs sont confrontés aux mêmes problèmes sémiotiques : *soutenir la conformité de la représentation qu'ils introduisent*. Tous ces mathématiciens introduisent des représentations nouvelles et soutenir la conformité de celles-ci les confronte à des problèmes qui ne sont peut-être pas reconnus comme des problèmes mathématiques mais qui doivent néanmoins être résolus et qui donnent à leurs textes certaines de leurs principales caractéristiques. Ces similitudes ne sont que les indices de l'*objectivité* de ces problèmes sémiotiques et d'une même *nécessité* de les résoudre. C'est ce qu'il s'agit peut-être surtout de reconnaître.

# II -Histoires inaugurales

Ce qui précède tente de montrer que la thèse et l'article de Turing partagent des caractéristiques remarquables avec un ensemble d'autres énoncés et de textes. Cette seconde partie va être consacrée à établir un second rapport entre ces énoncés et ces textes. Le premier rapport était de nature sémiotique et fonctionnel : la thèse de Turing partage avec d'autres énoncés le fait d'affirmer pour la première fois la conformité d'une représentation avec ce qu'elle représente et son article partage avec d'autres textes le fait d'être confronté aux problèmes sémiotiques spécifiques que posent le soutien d'une telle conformité. Le second rapport est de nature *historique*.

Toute représentation doit *nécessairement* avoir été inaugurée. Les textes inauguraux distingués ici ne correspondent qu'à un type particulier d'inauguration lié aux caractéristiques particulières de la représentation inaugurée. Ces caractéristiques n'impliquent pas seulement un type de texte : elles impliquent un processus historique qui s'étend au-delà d'eux. Cette conformité ne peut d'abord manquer d'être *contestée* : les textes inauguraux déterminent des *controverses spécifiques* relatives à la conformité, soit en contestant celle-ci soit au contraire en contestant l'impossibilité de la démontrer et en en proposant dès lors des démonstrations qui, immanquablement, pourront être contestées.

Après l'inauguration et les controverses qu'elle entraîne, un troisième temps peut encore être objectivement distingué[11] : celui de la reprise et de l'usage de la représentation introduite. Si la validité de l'énoncé inaugural ne saurait avoir été établie par les textes inauguraux, ils ont pu néanmoins permettre l'adoption de la représentation associée. L'enjeu et l'intérêt historiques de ces textes est l'introduction d'une nouvelle représentation au sein des mathématiques, c'est-à-dire la modification des moyens d'expression mathématiques. Ils sont un cadre privilégié pour reconnaître et étudier les questions qui s'y rapportent et déterminer les caractéristiques et les implications des représentations considérées. On peut en particulier distinguer les premiers textes qui *disposent* de la représentation inaugurée. Parfois ce sont les mêmes que ceux qui inaugurent la représentation, comme c'est le cas de Turing, mais pas toujours. Ces textes ont en commun de reprendre un système d'expression qui n'est pas encore généralement reçu. Il ne saurait être question de les caractériser comme il a été possible de le faire pour les textes inauguraux. Il s'agit ici surtout de reconnaître la différence entre *inaugurer* et *disposer* d'une telle représentation. Le statut de la représentation n'est pas le même dans les deux cas (même si ce changement peut être opéré dans un même texte). Cette différence peut *nécessiter* des transformations que les auteurs de ces textes doivent opérer comme les premiers avaient à les inaugurer. Ces transformations vont en particulier permettre des énoncés et des démonstrations pour lesquels l'inauguration ne suffit pas. Ce peut être aussi l'occasion de changer les dénominations utilisées dans le texte inaugural et de parler de « séries de Fourier » au lieu de « séries trigonométriques », de « machines de Turing » au lieu de « machines logiques », etc, changement qui ne peut être opéré que par un autre auteur et qui participe d'une autonomisation de la représentation par rapport à celui qui l'a inaugurée.

Il ne saurait être question de reprendre ici cette histoire même pour les seules machines logiques de Turing. Je me limiterai à nouveau à indiquer la similitude entre la réception des machines de Turing et celle d'autres représentations introduites dans d'autres textes inauguraux. Plus que cette similitude, c'est encore la raison de celle-ci qu'il importe de reconnaître. Mais comme elle tient toujours aux caractéristiques communes des systèmes d'expressions considérés, à laquelle la première partie a déjà été consacrée, je vais surtout m'attacher à montrer le *rapport historique* de l'énoncé et du texte inaugural de Turing avec d'autres inaugurations.

---

11 Il convient de distinguer aussi le temps qui précède l'inauguration de la représentation, c'est-à-dire tout le temps durant lequel une telle représentation n'était pas disponible, voire n'était souvent tout simplement pas concevable. Disposer de ces représentations a de nombreuses conséquences. Les ignorer conduit à autant d'anachronismes.

## 1 - La conformité controversée

La thèse de Turing, on le sait, a été et reste controversée. D'abord sceptique à l'égard de la formulation de Church, Gödel semble avoir été convaincu par la formulation de Turing (Gödel, postscriptum in Davis 1965, 71 ; Kleene 1981, 59 ; Davis 1982 ; Rosser 1984, 345 ; Kleene 1994, 44 ; Sieg 2006), sa position évoluant encore ensuite (Gödel 1972 ; Wang 1974 ; Webb 1990 ; Cassou-Noguès 2008). Cette thèse a aussi été contestée avec d'autres arguments à la fin des années 1950 par Laszlo Kalmar (1959), Rosza Péter (1959), J. Porte (1960), et leurs arguments contestés notamment par Mendelson (1963) et Moschovakis (1968). Elle a ensuite encore été contestée par Kreisel (1970, 1987) et Thomas (1973). Plus récemment encore de nouvelles critiques ont été développées à partir de la prise en compte des interactions entre calculateurs (Goldin & Wegner 2008).

La conformité de la représentation des problèmes de géométrie par les équations polynomiales n'a évidemment pas non plus manquée d'être contestée. James Gregory doutait que les moyens algébriques considérés soient suffisant pour la quadrature du cercle, problème géométrique par excellence :

> "J'ai parfois pensé, cher lecteur, que l'analyse était insuffisante avec ses cinq opérations et sa méthode générale, pour rechercher les proportions de toutes les quantités, comme on le voit affirmer par Descartes au début de sa Géométrie. Si, en effet, tel était le cas, on pourrait se vanter d'effectuer par ce moyen complètement la quadrature du cercle." Gregory 1667, cité et traduit in Dugac 2003, 42

Kepler (1619) avait déjà contesté la conformité de l'approche algébrique des problèmes géométriques (Bos 2001, 189-193). Cette conformité a aussi notoirement été contestée par Newton :

> «Les équations sont les expressions d'un calcul arithmétique et en tant que telles n'ont pas leur place en géométrie (...). Les multiplications, les divisions et autres opérations du genre ont été récemment introduites en géométrie, et cela sans précautions et contre les principes premiers de cette science (...). Ainsi, ces deux sciences [l'arithmétique et la géométrie] ne doivent pas être confondues. Les Anciens les tenaient distinctes avec tant d'attention qu'ils n'ont jamais introduit des termes arithmétiques en géométrie. Et les Modernes confondant l'une avec l'autre, ont perdu la simplicité qui fait toute l'élégance de la géométrie.» Newton, *Arithmetica Universalis*, appendice

La conformité ne requiert pas seulement la possibilité de mettre les problèmes en équations ; elle ne se réduit pas à la seule correspondance entre problèmes et équations. Elle suppose aussi que la résolution à partir des équations soit elle-même conforme à la construction géométrique du problème : les solutions obtenues algébriquement doivent coïncider avec les solutions obtenues géométriquement. La condition est double : toutes les solutions géométriques doivent être trouvées mais toutes les solutions trouvées doivent aussi être *a priori* reconnues comme étant géométriques. Or cela dépend des courbes que l'on s'autorise dans la résolution des problèmes : l'appréciation de cette conformité dépend de l'appréciation de la conformité de la représentation des courbes géométriques par les équations algébriques (deuxième énoncé inaugural de Descartes). Roberval a pu ainsi dénoncer le caractère incomplet de la solution proposée par Descartes au problème de Pappus (Maronne 2008). La conformité concerne aussi la simplicité : les solutions algébriquement simples doivent correspondre aux solutions géométriquement simples. Or la possibilité de rendre compte de la simplicité géométrique par le degré des équations a aussi été contestée, notamment par Newton et Jacques Bernoulli (Bos 1984, 358-366). Outre la simplicité, l'élégance devrait aussi être préservée, et Newton, à nouveau, conteste qu'elle le soit (Guicciardini 2002, 317).

L'hétérogénéité des termes mis en correspondance, l'absence d'une représentation uniforme de l'une des deux totalités en jeu et l'absence d'une représentation uniforme des propriétés visées par la conformité rendent en partie compte des controverses auxquelles les énoncés inauguraux de Turing et de Descartes ont donné lieu.

La conformité de le représentation des démonstrations et des propositions mathématiques soutenue dans les *Principia mathematica* de Whitehead & Russell n'a pas non plus manqué d'être contestée. Gödel après avoir tiré le meilleur parti de cette représentation des mathématiques avec ses théorèmes de complétudes et d'incomplétudes, comme avant lui Newton l'avait fait des représentations introduites par Descartes, n'en dénonça pas moins la conformité et s'efforça d'établir l'irréductibilité des mathématiques à une syntaxe (Wang 1981, 658 ; Gödel 1953 in Gödel 1995, 334 sq). Deux lettres adressées à Gödel montrent que Zermelo contestait lui aussi la conformité de cette représentation des mathématiques :

> « Exactement comme dans les paradoxes de Richard et de Skolem, l'erreur réside dans l'hypothèse (erronée) que toute notion mathématique définissable peut être exprimée par une « combinaison finie de signes » (selon un système *fixé*) – ce que j'appelle le «préjugé finitiste». En réalité, la situation est assez différente, et c'est seulement quand ce préjugé sera dépassé (tâche à laquelle je me suis attelé) qu'une « métamathématique » raisonnable sera possible. Interprétée de manière correcte, votre démonstration peut justement y contribuer beaucoup et rendre par conséquent un service substantiel à la cause de la vérité. » Zermelo à Gödel, 21 septembre 1931, cité in Dawson 1985, 69 (trad. AH)

Jean Dieudonné, s'exprimant pour Bourbaki, contesta aussi fermement la conformité de cette représentation des mathématiques :

> "Faire des mathématiques une partie de la logique (...), une telle affirmation est aussi absurde que celle qui consisterait à dire que les oeuvres de Shakespeare ou de Goethe font partie de la grammaire." Jean Dieudonné, "La philosophie des mathématiques de Bourbaki" in Dieudonné 1981, 31

L'approche intuitionniste des mathématiques ne peut aussi que contester la conformité de cette représentation des mathématiques en même temps que la possibilité d'être elle-même représentée de cette manière.

Sans doute convient-il de reconnaître les caractéristiques communes de ces controverses car elles doivent se retrouver dans les analyses qui en sont faites, même s'il ne saurait être question de les y réduire.

## 2 - Disposer et mettre à disposition

Tous les textes et les énoncés inauguraux soutiennent la conformité de la représentation qu'ils introduisent et celle-ci est comme on vient de le voir inévitablement contestée. Cela n'empêche pas ces représentations d'être ensuite reprises, en étant ou non tenues pour conformes.
Les déclarations en faveur de la thèse de Turing ne manquent pas. Citons simplement :

> "De nos jours, la thèse de Church peut difficilement être contestée" Davis 1982, 22 (trad. AH)

De manière semblable, la conformité de la représentation algébrique des problèmes de géométrie a eu ses adeptes comme Pierre Hérigone (1634) ou encore Jacques Ozanam (1691)[12], ce dernier commençant par exemple son article « Algèbre » en affirmant :

> « L'Algèbre est une science, par le moyen de laquelle on peut resoudre tout Probleme possible dans la Mathematiques. » Ozanam 1691, 61

Mais surtout, les différentes représentations des nombres et des fonctions calculables (machines de Turing, lambda-définissabilité, fonctions générales récursives, etc.) ont pu être utilisées pour leur intérêt propre sans que leur conformité ne soit plus un préalable. Rogers, auteur de l'un des ouvrages de référence sur la théorie de la récursivité, peut ainsi déclarer :

> « Il doit cependant être souligné que notre théorie tire sa signification principale de son intérêt pour les mathématiques pures. Elle fournit des structures d'une beauté naturelle et intrinsèque remarquable. Elle ouvre des perspectives nouvelles et souvent profondes sur d'autres parties. Ces perspectives ont été particulièrement utiles en logique mathématique, et elles ont de la même manière été de plus en plus utiles dans des parties plus classiques. » Rogers 1967, 31 (trad. AH)

Ces différentes représentations ont aussi systématiquement donné lieu à diverses extensions ou généralisations qui supposent l'abandon du souci de leur conformité.
Les généralisations des machines de Turing ne manquent pas : automates cellulaires, machines à signaux, etc. Les expressions polynomiales ont de la même manière été étendues aux séries et celles-ci aux séries trigonométriques. Et de la même manière, après que Fourier ait soutenu la possibilité d'exprimer par de telles séries toutes les fonctions, la théorie des séries s'est développée et a été généralisée (théorie des distributions, ondelettes) en abandonnant et en oubliant ce souci de conformité, voir en oubliant leurs interprétations physiques :

> « J'ai travaillé pendant de nombreuses années sur la théorie des intégrales de Fourier sous la direction [de Hardy] avant de découvrir que cette théorie avait des applications en mathématiques appliquées, pour autant que la solution de certaines équations différentielles puisse être qualifiée d' « appliquée ». » Titchmarsh 1950, 85 (trad. AH)

La représentation des propositions et des démonstrations inaugurée dans les *Principia mathematica* a comme on l'a vu aussi été reprise, avec d'autres, et a été pour cela

---

12 Dans les deux cas il s'agit du système d'expression inauguré par Viète.

*transformée*. L'article de Post (1921) en est un exemple. Il ne considère en l'occurrence qu'une partie de ce système, le calcul propositionnel, dont il énonce et démontre la complétude. Un tel théorème n'a pas sa place dans les *Principia mathematica* : c'est un théorème *sur* le système d'expression utilisé et non *dans* celui-ci. Pour les auteurs des *Principia* la démonstration d'un tel résultat n'est ni utile ni possible. En sortant de leur cadre, une telle démonstration est incompatible avec leur intention d'établir que toutes les propositions mathématiques peuvent être exprimées *dans* ce système. En outre, si la démonstration de Post porte bien *sur* cette représentation, elle n'est pas et ne saurait être exprimée *dans* celle-ci. Ni son énoncé, ni sa démonstration ne sauraient donc être intégrés aux *Principia*. Mais pour énoncer et démontrer son théorème Post a besoin de donner une présentation d'ensemble du calcul propositionnel, ce que les auteurs des *Principia* n'avaient pas fait. Une telle présentation *nécessaire* pour traiter de la complétude *nécessite* de changer le statut de la représentation inaugurée dans les *Principia*.

Jacques Herbrand (1930, oeuvres, 35) tient aussi pour établie la conformité de la représentation inaugurée par Whitehead & Russell. Comme Post, il ne s'agit plus pour lui de soutenir la conformité de la représentation inaugurée par Whitehead & Russell ni d'ajouter des propositions supplémentaires à celles des *Principia mathematica* mais de proposer de nouveaux développements *à partir de* cette représentation :

> « Arrivés à ce point où l'on découvre un parallélisme absolu entre le raisonnement mathématique et les combinaisons des signes, il devient naturel d'étudier ce système de signes pour lui-même et de se poser des problèmes à son sujet; et leur solution aura immédiatement un écho dans notre connaissance générale des mathématiques (c'est ainsi que l'on peut démontrer la non-contradiction de certaines théories). Nous nous poserons des problèmes du genre suivant : les propositions possèdent-elles telle propriété ? Ou : les propositions vraies possèdent-elles telle propriété ? Nous pourrons également chercher un critère permettant de reconnaître sûrement si une proposition donnée est vraie dans une théorie déterminée ; sa découverte dans un cas suffisamment général serait évidemment d'une importance théorique considérable. » Herbrand 1930 *in* Herbrand 1968, 37

Jacques Herbrand remarque lui-même que les démonstrations de ses théorèmes recourent soit à des *récurrences* sur la construction des propositions soit à des *récurrences* sur leur démonstration[13] : les théorèmes et les démonstrations qu'il propose ne dépendent pas de la conformité mais des caractéristiques de la représentation qu'elle a permise d'introduire. En l'occurrence, la représentation introduite dans les *Principia mathematica* met à la disposition de Herbrand un système d'expressions qui rend possible des démonstrations *par récurrence* auparavant impossibles. Les caractéristiques de la représentation inaugurée rendent cette fois possibles de nouveaux développements mathématiques. Mais il convient aussi d'observer la transformation qui s'opère entre les *Principia mathematica* et les développements proposés par le jeune logicien. Celui-ci peut en effet considérer comme une totalité constituée les propositions définies dans les *Principa mathematica* alors que leurs auteurs avaient eux à faire la preuve de la conformité de leur représentation ; ils ne pouvaient en l'occurrence à la fois administrer cette preuve et tirer parti de leur représentation. Les propositions et les démonstrations des *Principia* ont permis d'inaugurer des représentations dont Herbrand a pu ensuite disposer. Cette différence et le lien entre les caractéristiques des représentations utilisées et les démonstrations proposées ressortent ici particulièrement bien : les démonstrations des *Principia mathematica* n'utilisent que la substitution et le *modus ponens,* celles d'Herbrand des récurrences (sur les représentations adoptées). Les démonstrations dépendent des représentations dont on dispose et changent avec elles.

---

[13] Il distingue un troisième type de démonstrations qui n'en est qu'une généralisation mais dans les faits fondés sur les deux premiers relatif à « des séries de propositions dont chacune se déduit de la précédente par une opération déterminée ; ce qui fait qu'on peut numéroter ces propositions avec des chiffres » Herbrand 1930 in Herbrand 1968, 39.

Les articles dans lesquels Gödel démontre ses théorèmes de complétude et d'incomplétude sont bien sûr aussi des exemples à la fois de reprise et de transformation de la représentation inaugurée dans les *Principia mathematica.* Mais pour lui, cette représentation ne formalise pas toutes mais seulement de « larges secteurs » des mathématiques :

> « Le développement des mathématiques vers plus d'exactitude a conduit, comme nous le savons, à en formaliser de larges secteurs, de telle sorte que la démonstration puisse s'y effectuer uniquement au moyen de quelques règles mécaniques. Les systèmes formels les plus complets établis jusqu'à ce jour sont, d'un côté, le système des P*rincipia Mathematica* (PM) et, de l'autre, le système axiomatique de la théorie des ensembles établi par Zermelo-Fraenkel (et développé par J. von Neumann). » Gödel 1931, trad. 107

Tous ces travaux reprennent la représentation inaugurée par Whitehead & Russell. Ils ne soutiennent pas nécessairement la conformité qui en était la raison d'être. Celle-ci entre bien néanmoins toujours pour une part dans l'attrait de cette représentation. La portée couramment donnée aux théorèmes d'incomplétude de Gödel l'atteste.

L'article de Turing est lui-même en partie fondé sur le fait de disposer de représentations conformes des propositions et des démonstrations mathématiques[14]. Il s'inscrit lui aussi dans l'histoire de la réception de ces représentations.

La plupart de ses arguments pour soutenir la conformité de sa représentation des nombres calculables font appel à la conformité de la représentation des propositions et des démonstrations, en l'occurrence celle donnée par Hilbert & Bernays (1934). C'est là le second rapport, de nature historique cette fois, entre les inaugurations. Examinons-le plus précisément.

Turing soutient la conformité de sa représentation au moyen de quatre arguments (qu'il répartit suivant trois types) :

1. les machines logiques reproduisent les opérations effectuées par un humain qui calcule ;
2. l'équivalence avec une autre définition de la calculabilité ;
3. les machines logiques procèdent comme un calculateur dissipé ;
4. le calcul de classes de nombres calculables.

La conformité de la représentation des propositions et des démonstrations intervient dans trois de ces quatre arguments (2, 3 et 4).

Le système d'expression de Hilbert permet à Turing de donner une définition de la calculabilité différente mais équivalente à la sienne « au cas où la nouvelle définition [celle à partir du système d'Hilbert!] aurait un sens intuitif plus évident ». Cette définition, déjà évoquée, consiste à considérer qu'un nombre α, compris entre 0 et 1 et dont le développement décimal est exprimé en binaire, est calculable si pour tout entier *n* il est possible de démontrer dans le système de Hilbert soit la formule ($A_n$) disant que la $n^{\text{ème}}$ décimale de α est 1 soit la formule ($B_n$) disant que la $n^{\text{ème}}$ décimale de α est 0. A partir d'une machine logique qui énumère toutes les formules démontrables du système de Hilbert, Turing peut décrire une machine qui trouve toutes les formules des nombres calculables (suivant cette seconde définition), ce qui montre que les nombres calculables (toujours suivant cette définition) sont calculables au sens de Turing. L'idée consiste à rendre la possibilité ou non de calculer un nombre par la possibilité ou non de démontrer la propositions disant que la $n^{\text{ème}}$ décimale vaut 1 ou celle disant qu'elle vaut 0. Même si Turing ne fait qu'en esquisser la démonstration, cette équivalence est tout à fait susceptible d'une démonstration complète. Là n'est pas la question. En revanche, il ne s'agit pas seulement de démontrer cette équivalence, mais de proposer un argument pour soutenir la conformité de la représentation proposée de la calculabilité. Cette équivalence n'offre un tel argument que si le système de Hilbert donne elle-même une expression conforme d'une part de toutes les manières (formules) de définir les décimales d'un nombre et d'autre part de toutes les démonstrations susceptibles de les démontrer. Dans les deux cas, cela suppose que ce système d'expressions donne toutes et exactement toutes les propositions ou démonstrations, sans en ajouter aucune. C'est donc en en appelant à la conviction que peut avoir son lecteur de la conformité de la représentation proposée par Hilbert que Turing soutient la conformité de sa représentation de la calculabilité.

---

14 Ces liens se retrouvent aussi dans les travaux de Church et de Post.

Considérons à présent son deuxième argument, celui du calculateur dissipé, qui fait à nouveau intervenir la conformité de la représentation cette fois uniquement des propositions. Il s'agit, comme pour le premier qui analysait les opérations effectuées par un calculateur humain, de montrer que ses machines permettent de représenter une situation à laquelle tout calcul est susceptible d'être ramené. Il imagine pour cela une personne qui fait un calcul en étant systématiquement interrompue à chaque opération. Pour mener à bien son calcul, cette personne doit disposer d'une note contenant toutes les informations nécessaires pour reprendre le calcul interrompu. Il suffit ensuite d'admettre :

> « que cette dépendance peut s'exprimer dans le formalisme du calcul fonctionnel. Autrement dit, nous admettons qu'il existe un axiome $\mathcal{U}$ qui exprime les règles que doit suivre le calculateur en face d'une formule d'état de n'importe quelle étape en vue d'obtenir la suivante. » Turing 1936, trad. fr. 84

Comme on le voit, Turing admet ici que cette note pourra toujours être écrite dans le système de Hilbert.

Considérons maintenant le troisième argument recourant à la conformité de cette représentation. Pour établir que des classes entières de nombres sont calculables, et non plus simplement certains nombres particuliers, Turing démontre quelques théorèmes généraux élémentaires sur les *fonctions* calculables : la somme de deux fonctions calculables est calculable, la composée d'une fonction calculable par une fonction calculable est calculable, etc. La démonstration de ces théorèmes a besoin de la représentation d'une fonction calculable *quelconque* $\phi(x,y)=z$ (Turing 1936, trad. fr. p. 88 sq). Turing va à nouveau considérer ici qu'il existe pour n'importe quelle fonction calculable $\phi$ une formule du système d'expressions logique d'Hilbert exprimant l'assertion $\phi(x,y)=z$. La démonstration du théorème est ainsi ramenée à la possibilité de démontrer une formule donnée dans le système de Hilbert (la possibilité de cette démonstration rendant toujours compte de la calculabilité). Et à nouveau Turing en appelle ici à la conviction qu'il dispose d'un système d'expression susceptible d'exprimer *n'importe quelle* relation $\phi(x,y)=z$.

Le fait de disposer d'une représentation conforme des propositions et des démonstrations mathématiques est ainsi un des principaux ressorts de la plupart des arguments donnés par Turing pour soutenir la conformité de sa représentation des nombres calculables[15]. L'antériorité de cette inauguration n'est donc pas contingente. La concomitance d'articles soutenant une « thèse » équivalente ne l'est sans doute pas non plus. Les énoncés et les textes inauguraux ne font pas que présenter des caractéristiques communes : ils ont entre eux des rapports de présupposition objectifs. Il est ainsi possible de déterminer les représentations conformes dont ils *disposent.*

La représentation des démonstrations et propositions mathématiques n'est pas la seule représentation conforme dont dispose Turing. Il dispose aussi de celle des nombres réels compris entre 0 et 1. Elle peut aujourd'hui sembler aller de soi, mais elle n'est pourtant peut-être pas moins remarquable. Sans elle, le problème d'une représentation conforme de la calculabilité ne se poserait pas[16]. Sans elle, il ne pourrait tout simplement pas y avoir de machine calculant la somme de deux nombres calculables. En effet, une telle machine requiert une représentation de *tous* les nombres par un système d'expressions réduit à un nombre fini et fixé à l'avance d'expressions (les chiffres) : il faut connaître *a priori* les symboles nécessaires à l'écriture de ces nombres pour pouvoir écrire toutes les *m*-configurations nécessaires à la description complète d'une machine s'appliquant à *tous* les nombres. Sans cela, l'addition de deux nombres nécessiterait de faire appel à un oracle! Cette condition réalisée par le système de numération de position utilisé par Turing ne l'est pas, loin s'en faut, par tous les systèmes de numération qui ont généralement besoin de *nouveaux* symboles pour exprimer des nombres aussi grands ou proches de zéro que l'on veut. Autrement dit, les machines logiques de Turing représentent les nombres calculables à condition de les exprimer dans un système de numération ayant certaines caractéristiques elles-mêmes tout à fait particulières et remarquables. C'est aussi en raison du choix de ce système d'expression que les machines de Turing ont des bandes indéfiniment prolongeables à droite. La thèse de Turing doit donc aussi être restituée dans l'histoire de la réception de ce système de numération dont elle dispose et qu'elle présuppose. Son article est d'ailleurs émaillé de considérations sur le système de numération positionnel, dont dépend en fait son analyse des opérations effectuées par un calculateur. Il ajoutera d'ailleurs un correctif pour tenir compte du problème posé par l'existence pour certains nombres de plusieurs expressions différentes dans ce système, ce qui montre la dépendance de ses machines à ce *système de représentation* des nombres.

Il est ainsi illusoire de croire avec Gödel que la définition de ces machines soit indépendante du formalisme choisi[17] : son adéquation avec le formalisme utilisé pour

---

[15] L'*Entscheidungsproblem* que résout Turing est aussi un problème qui n'est envisageable que si l'on dispose d'une représentation de toutes les démonstrations et propositions mathématiques.

[16] Sur l'importance de la représentation des nombres dans la mécanisation des calculs voir Bullynck à paraître, 41.

[17] "Tarski a souligné dans sa conférence (et je crois à juste titre) l'extrême importance du concept de récursivité générale (ou de calculabilité au sens de Turing). Il me semble que cette importance tient

représenter les nombres permet seulement d'en masquer la dépendance. L'indépendance n'est ici qu'une dépendance ignorée.

La *traduction* française de l'article de Turing offre un exemple anecdotique mais néanmoins significatif du statut de la représentation logique des propositions mathématiques.
A gauche ci-dessous un extrait de l'article de Turing, à droit sa traduction française.

| | |
|---|---|
| We shall say that a sequence $\beta_n$ of computable numbers *converges computably* if there is a computable integral valued function $N(\varepsilon)$ of the computable variable $\varepsilon$, such that we can show that, if $\varepsilon > 0$ and $n > N(\varepsilon)$ and $m > N(\varepsilon)$ then $\|\beta_n - \beta_m\| < \varepsilon$. | Nous dirons qu'une séquence $\beta_n$ de nombres calculables *converge pour le calcul* s'il existe une fonction calculable à valeur entière $N(\varepsilon)$ de la variable calculable $\varepsilon$, telle que : $\forall \varepsilon > 0, \forall m > N(\varepsilon), \forall n > N(\varepsilon), \|\beta_n - \beta_m\| < \varepsilon$ |
| Turing 1936, 256 | Turing 1936, trad. fr. 87 |

La condition que Turing exprime par « if $\varepsilon > 0$ and $n > N(\varepsilon)$ and $m > N(\varepsilon)$ then $\|\beta_n - \beta_m\| < \varepsilon$ » a été dans la traduction isolée par un alinéa du texte et surtout, elle a été traduite par « $\forall \varepsilon > 0, \forall m > N(\varepsilon), \forall n > N(\varepsilon), \|\beta_n - \beta_m\| < \varepsilon$ ». Le traducteur a exprimé de manière entièrement symbolique ce que Turing exprimait dans un mélange de langue naturelle et de symbolisme. La notation de la quantification universelle utilisée par le traducteur n'est pas celle qu'utilise Turing. Mais surtout, Turing n'utilise pas cette représentation si elle n'est pas *nécessaire* à son propos : il n'y recourt pas dans le seul but, comme le fait ici le traducteur, d'abréger l'expression d'une condition. La *fonction* de la représentation des propositions logiques n'est pas la même pour lui et pour le traducteur. Ce n'est pas encore pour Turing le mode d'expression par défaut ou un moyen commode d'exprimer une proposition : *la différence entre les propositions et leur représentation symbolique est chez lui encore bien marquée*. L'assimilation de la représentation à ce qu'elle représente est plus avancée chez le traducteur. Cette représentation des propositions logiques a acquis pour lui un statut semblable à bien d'autres représentations mathématiques, ce qui lui permet de la mêler à d'autres, ce qui n'était pas encore le cas pour Turing. On a pourtant vu qu'il disposait aussi de cette représentation et que sa conformité jouait un rôle essentiel dans son article. Mais il n'en dispose pas comme le traducteur. En particulier, il n'isole pas, comme celui-ci le fait ici avec le symbole de quantification, une expression du *système d'expressions* auquel elle est encore pour lui attachée.
Les *Principia mathematica*, cet article de Turing et sa traduction marquent trois moments différents de l'histoire de la réception de cette représentation. Elle a été introduite en tant que représentation conforme, reprise ou non en tant que telle, et elle finie par être à peu près complètement transparente, comme le système de numération positionnel l'est dans

---

largement au fait qu'avec ce concept nous avons réussi pour la première fois à donner une définition absolue d'une notion épistémologique intéressante, *i.e.* qui ne dépende pas du formalisme choisi. Dans tous les autres cas considérés avant, comme la démontrabilité ou la définissabilité, nous n'avons été capables de les définir que relativement à un langage donné, et il est clair que pour chaque langage individuel celle obtenue n'est pas celle que nous cherchions. En revanche, pour le concept de calculabilité, bien qu'il ne s'agisse essentiellement que d'un cas particulier de démontrabilité ou de décidabilité la situation est différente. Par une sorte de miracle il n'est pas nécessaire de distinguer des ordres, et la procédure diagonale ne nous fait pas sortir de la notion définie." Gödel, Kurt, "Remarks before the Princeton bicentennial conference on problems in mathematics", 1946, cité *in* Gandy 1988, 80 (trad. AH) l'article de Turing[18].

# Conclusion

   La définition qui a été proposée de la notion d'*énoncé inaugural* distingue un ensemble de caractéristiques à l'oeuvre dans l'énoncé de la thèse de Turing. Ces caractéristiques s'apprécient ensemble, la valeur de chacune venant en partie de la satisfaction des autres. La conformité de la représentation introduite, c'est-à-dire le fait qu'elle permette d'exprimer toutes les propriétés de ce qu'elle représente, est néanmoins sans doute la plus exorbitante et celle qui confère à ces énoncés leur caractère exceptionnel. Mais si la thèse de Turing est bien de ce point de vue aussi un énoncé exceptionnel, il n'est pas unique. Des inaugurations de représentations conformes jalonnent l'histoire des mathématiques comme la représentation des problèmes et des courbes géométriques par des équations, la représentation des fonctions par des séries trigonométriques, la représentation des propositions et des démonstrations mathématiques par des formules logiques, la représentation des points d'une droite par les nombres réels, la représentation des ensembles de nombres définissables par des formules ou même simplement la représentation des nombres entiers puis réels dans le système de numération décimal de position.. La liste n'est évidemment pas close  (il faut au moins mentionner ici la représentation conforme des objets mathématiques par des ensembles), mais il devrait néanmoins apparaître à chacun qu'une telle liste peut être établie complètement. Ces énoncés marquent tous des temps forts et déjà, indépendamment, reconnus comme tels de l'histoire des mathématiques. Reconnaître que la thèse de Turing est un énoncé inaugural *parmi d'autres* permet néanmoins d'inscrire la représentation des nombres calculables par les machines de Turing dans la liste remarquable des représentations conformes qui ont à chaque fois marqué un moment important de l'histoire des mathématiques. Il ne s'agit pas seulement d'une liste plus ou moins arbitraire d'énoncés semblables : les caractéristiques des représentations qu'ils inaugurent *rendent compte* de leur importance. Elles conduisent aussi à distinguer au côté de la notion d'*énoncé* inaugural celle de *texte* inaugural : Turing n'introduit pas seulement un type d'énoncé remarquable, il écrit aussi un texte dont les caractéristiques sont *déterminées* par celles de la représentation qu'il inaugure. La conformité des représentations inaugurées détermine un type particulier d'énoncés mais aussi de textes.
L'énoncé et l'article de Turing ne sont pas seulement des exemples d'énoncés et de textes inauguraux. Des liens historiques existent entre les représentations inaugurées :  Turing tient pour acquis que lui et son lecteur disposent d'une représentation conforme des démonstrations et des propositions mathématiques mais aussi des nombres réels. Son article s'inscrit donc doublement dans l'histoire de ces inaugurations : c'est un texte inaugural qui dispose de représentations déjà inaugurées et dont il peut en l'occurrence continuer de faire valoir la conformité.

---

18 Des évolutions semblables peuvent être mises en évidence à partir des *traductions* des textes qui inaugurent la représentation des nombres par le système décimal de position, voir Herreman (2001).

# Bibliographie